\documentstyle[amscd,amssymb,verbatim,12pt]{amsart}
\newcommand{\Aff}{\mbox{\rm Aff}}
\newcommand{\cen}{\mbox{\rm center}}
\newcommand{\diam}{\mbox{\rm diam}}
\newcommand{\Div}{\mbox{\rm div}}

\newcommand{\Ker}{\mbox{\rm Ker}}

\newcommand{\Image}{\mbox{\rm Im}}
\newcommand{\Aut}{\mbox{\rm Aut}}

\newcommand{\End}{\mbox{\rm End}}

\newcommand{\ind}{\mbox{\rm ind}}

\newcommand{\Ric}{\mbox{\rm Ric}}

\newcommand{\vol}{\mbox{\rm vol}}

\newcommand{\Isom}{\mbox{\rm Isom}}
\newcommand{\isom}{\mbox{\rm isom}}

\newcommand{\Diff}{\mbox{\rm Diff}}

\newcommand{\R}{{\Bbb R}}
\newcommand{\C}{{\Bbb C}}

\newcommand{\Z}{{\Bbb Z}}

\newcommand{\hA}{\widehat{A}}

\theoremstyle{plain}
\newtheorem{definition}{Definition}

\newtheorem{lemma}{Lemma}

\newtheorem{proposition}{Proposition}

\newtheorem{question}{Question}
\numberwithin{equation}{section}
\renewcommand{\rm}{\normalshape}
\begin{document}
\title{$\hA$-genus and Collapsing}
\author{John Lott}
\address{Department of Mathematics\\
University of Michigan\\
Ann Arbor, MI  48109-1109\\
USA}
\email{lott@@math.lsa.umich.edu}
\thanks{Research supported by NSF grant DMS-9403652}
\date{April 24, 1997}
\maketitle
\begin{abstract}
If $M$ is a compact spin manifold, we give relationships between the
vanishing of $\hA(M)$ and the possibility that $M$ can collapse with
curvature bounded below.
\end{abstract}
\section{Introduction}
The purpose of this paper is to extend the following simple lemma.
\begin{lemma} \label{lemma1}
If $M$ is a connected
closed Riemannian spin manifold of nonnegative sectional curvature with
$\dim(M) > 0$ then $\widehat{A}(M) = 0$.
\end{lemma}
\begin{pf}
Let $K$ denote the sectional curvature of $M$ and let $R$ denote its scalar 
curvature.
Suppose that $\widehat{A}(M) \ne 0$. Let $D$ denote the Dirac operator on $M$.
From the Atiyah-Singer index theorem, there is a nonzero 
spinor field $\psi$ on $M$ such that $D \psi = 0$. 
From Lichnerowicz' theorem,
\begin{equation}
0 = \int_M | D\psi|^2 \: d\vol =
\int_M | \nabla \psi |^2 \: d\vol + 
\int_M \frac{R}{4} \: |\psi|^2 \: d\vol.
\end{equation}
From our assumptions, $R \ge 0$. Hence $\nabla \psi = 0$. This implies that
$|\psi|^2$ is a nonzero constant function on $M$ and so we must also have
$R = 0$. Then as $K \ge 0$, we must have $K = 0$. This implies, from
the integral formula for $\widehat{A}(M)$
\cite[p. 231]{Lawson-Michelsohn (1989)}, that $\widehat{A}(M) = 0$. 
\end{pf}

The spin condition is necessary in Lemma \ref{lemma1}, as can be seen
in the case of $M = \C P^{2k}$.

\begin{definition}
A connected closed manifold $M$ is almost-nonnegatively-curved if for every
$\epsilon > 0$, there is a Riemannian metric $g$ on $M$ such that
$K(M, g) \cdot \diam(M, g)^2 \ge - \epsilon$.
\end{definition}

Special examples of almost-nonnegatively-curved manifolds are given by
almost-flat manifolds; these all have vanishing $\hA$-genus, as can be
seen by the integral formula. Along with Lemma \ref{lemma1}, this raises the
following question.

\begin{question} \label{ques1}
Given $n \in \Z^+$, is there an $\epsilon(n) > 0$ such that if $M$ is a
connected closed Riemannian spin manifold with $K(M,g) \cdot \diam(M,g)^2
\ge - \epsilon(n)$ then $\hA(M) = 0$? 
\end{question}

We answer Question \ref{ques1} under the assumption of an upper
curvature bound.

\begin{proposition} \label{prop1}
For any $n \in \Z^+$ and any $\Lambda > 0$, there is an
$\epsilon(n, \Lambda) > 0$ such that if $M$ is a connected closed
$n$-dimensional Riemannian spin manifold with
\begin{equation}
-\epsilon(n, \Lambda) \le K(M,g) \cdot \diam(M,g)^2 \le \Lambda
\end{equation}
then $\hA(M) = 0$.
\end{proposition}

The proof of Proposition \ref{prop1} uses Gromov's convergence theorem
\cite{Gromov-Lafontaine-Pansu (1978),Kasue (1989)}.
Using the results of \cite{Petersen-Wei-Ye (1995)}, the upper bound on the
sectional curvature in Proposition \ref{prop1} can be replaced by
a lower bound on the conjugacy radius.

An affirmative answer to Question \ref{ques1} would imply that an
almost-nonnegatively-curved spin manifold has vanishing $\hA$-genus.
There is a fiber bundle construction to create new 
almost-nonnegatively-curved manifolds out of old ones.  The following
simple proposition shows that
the vanishing of the $\hA$-genus is consistent with this construction.
\begin{proposition} \label{prop2}
Let $N$ be a connected closed manifold of nonnegative sectional curvature.
Let $G$ be a compact Lie group which acts on $N$ by isometries.
Let $P$ be a principal $G$-bundle
with connected closed base $B$.  Put $M = P \times_G N$.\\
1. If $B$ is almost-nonnegatively-curved then
$M$ is almost-nonnegatively-curved
\cite[Theorem 0.18]{Fukaya-Yamaguchi (1992)}.\\
2. If $M$ is spin and $\dim(N) > 0$ then $\widehat{A}(M) = 0$.
\end{proposition}

Proposition \ref{prop2} covers a wide range of almost-nonnegatively-curved 
manifolds.
It seems conceivable that every almost-nonnegatively-curved manifold has
a finite cover which is the total space of a fiber bundle whose
base is almost-flat and whose fiber has a metric of nonnegative sectional
curvature.

Rescaling metrics, a manifold $M$ is 
almost-nonnegatively-curved if for every $\epsilon > 0$, there is
a Riemannian metric $g$ on $M$ such that $K(M, g) \ge -1$ and 
$\diam(M, g) \le \epsilon$. That is, there is a sequence of metrics
$\{g_i\}_{i=1}^\infty$ such that $K(M, g_i) \ge -1$ and the metric spaces
$\{(M, g_i)\}_{i=1}^\infty$ converge in the Gromov-Hausdorff topology to
a point.  It is natural to extend Question \ref{ques1} to a question about
the $\hA$-genus of a spin manifold $M$ with a sequence of metrics
$\{g_i\}_{i=1}^\infty$ such that $K(M, g_i) \ge -1$ and the metric spaces
$\{(M, g_i)\}_{i=1}^\infty$ converge in the Gromov-Hausdorff topology to
some lower-dimensional length space, not necessarily a point. 
The following definition is convenient for our purposes.

\begin{definition}
A connected manifold $M$ collapses with curvature bounded below
and diameter bounded above if there is a number $D > 0$ such that for any
$\epsilon > 0$, there is a 
Riemannian metric $g$ on $M$ with $K(M, g) \ge -1$,
$\diam(M, g) \le D$ and $\vol(M, g) \le \epsilon$.
\end{definition}

We remark that in the noncollapsing case there is a finiteness result
\cite{Grove-Petersen-Wu (1990)}.
Namely, given $D, v > 0$ and $n > 3$, there is a finite number of
homeomorphism classes of connected
manifolds $M^n$ admitting a Riemannian metric $g$
satisfying $K(M, g) \ge -1$, $\diam(M, g) \le D$ and $\vol(M, g) \ge v$.

\begin{question} \label{ques2}
Given $n \in \Z^+$ and $D > 0$, is there a $v(n, D) > 0$ such that
if $M$ is a connected closed $n$-dimensional Riemannian spin manifold with 
$K(M,g) \ge -1$, $\diam(M,g) \le D$ and $\vol(M,g) \le v(n, D)$ then
$\hA(M) = 0$?
\end{question}

An affirmative answer to Question \ref{ques2} would imply that a
spin manifold which collapses with curvature bounded below and
diameter bounded above has vanishing $\hA$-genus. In the next proposition
we show that this is indeed the case for a large
class of collapsing examples.

\begin{proposition} \label{prop3}
Let $Z$ and $N$ be connected closed Riemannian manifolds.  
Suppose that $N$ has nonnegative sectional curvature.
Let $G$ be a compact Lie group which acts on $Z$ and $N$ by isometries.
Suppose that for
a generic point $z$ in $Z$, the stabilizer group $G_z$ does not act
transitively on $N$. Suppose that the diagonal action of $G$ on
$Z \times N$ has the property that all of its orbits are principal orbits.
Let $M$ be the quotient manifold $Z \times_G N$. Then\\
1. $M$ collapses with curvature bounded below and diameter bounded above.
The collapsing sequence converges in the Gromov-Hausdorff topology to 
the length space $Z/G$.\\
2. If $M$ is spin then $\widehat{A}(M) = 0$.
\end{proposition}

Some special cases of Proposition \ref{prop3} are :\\
1. If $G$ acts freely on $Z$. Then Proposition \ref{prop3} is equivalent
to Proposition \ref{prop2}.\\
2. If $N = G$ is a connected compact Lie group which acts nontrivially on $Z$.
Then $M = Z \times_G G = Z$
and the second part of Proposition \ref{prop3} is equivalent to the 
Atiyah-Hirzebruch theorem \cite{Atiyah-Hirzebruch (1970)}. 

To put the results of this paper in perspective, let us mention
known necessary conditions 
for a connected closed manifold $M$ to be almost-nonnegatively-curved :\\
1. The fundamental group $\pi_1(M)$ must be virtually nilpotent
\cite[Theorem 0.1]{Fukaya-Yamaguchi (1992)}.\\
2. If $\pi_1(M)$ is infinite then the Euler characteristic of $M$ must
vanish \cite[Corollary 0.12]{Fukaya-Yamaguchi (1992)}.\\
3. $M$ must be dominated by a CW-complex with the number of cells bounded
above by a function of $\dim(M)$ \cite{Gromov (1981), Weiss (1996)}.\\
4. If $M$ is spin
then $|\hA(M)| \le 2^{\frac{dim(M)}{2} - 1}$. (This is a necessary
condition for $M$ to have almost-nonnegative-Ricci curvature 
\cite{Gallot (1983)}.)\\

I thank Peter Petersen for his interest in these questions.

\section{Proof of Proposition \ref{prop1}} 
For background material on spin geometry, we refer to
\cite{Lawson-Michelsohn (1989)}. 
Before giving the proof of Proposition \ref{prop1}, we must
discuss how to compare spinors on diffeomorphic Riemannian manifolds which
are not necessarily isometric. This is an elementary point which has caused 
confusion in the literature.

Let $M$ be a smooth connected closed $n$-dimensional oriented manifold. 
Let $PM$ be a principal $Spin(n)$-bundle on $M$.
Let $S_n$ be the complex spinor module of
$Spin(n)$. Then we can form the associated Hermitian
vector bundle $S = PM \times_{Spin(n)} S_n$ on $M$. The corresponding
spinor fields
are defined to be the sections of $S$, or equivalently, the
$Spin(n)$-equivariant maps from $PM$ to $S_n$. So far we have made no 
reference to a Riemannian metric on $M$.

Let $p : FM \rightarrow M$ be the oriented frame bundle of $M$, 
a principal $GL^+(n, \R)$-bundle on $M$.
Given $\gamma \in GL^+(n, \R)$, let $R_\gamma \in \Diff(FM)$ denote the right 
action of $\gamma$ on $FM$. There is a canonical $\R^n$-valued $1$-form 
$\theta$ on $FM$ such that if $f = \{f_i\}_{i=1}^n$ is an oriented frame at
$m \in M$ and $v \in T_fFM$ then $dp(v) = \sum_{i=1}^n \theta^i(v) f_i$.
It has the properties that\\
1. If $V$ is a vertical vector field on $FM$ then $\theta(V) = 0$.\\
2. For all $\gamma \in GL^+(n, \R)$, 
$R_\gamma^* \theta = \gamma^{-1} \cdot \theta$.\\
3. For all $f \in FM$, $\theta : T_fFM \rightarrow \R^n$ is onto. \\

Giving a Riemannian metric $g$ on $M$ is equivalent to giving a reduction 
$i : OM \rightarrow FM$ of the oriented frame bundle from a principal
$GL^+(n, \R)$-bundle to
a principal $SO(n)$-bundle $OM$. As a topological fiber bundle, $OM$ is
unique. We obtain an $\R^n$-valued $1$-form $\tau =
i^* \theta$ on $OM$ with the properties that\\
1. If $V$ is a vertical vector field on $OM$ then $\tau(V) = 0$.\\
2. For all $\gamma \in SO(n)$, 
$R_\gamma^* \tau = \gamma^{-1} \cdot \tau$.\\
3. For all $f \in OM$, $\tau : T_fOM \rightarrow \R^n$ is onto. \\

Conversely, given the topological $SO(n)$-bundle $\pi : OM \rightarrow M$ 
and an $\R^n$-valued $1$-form $\tau$ on $OM$ satisfying properties 1.-3. 
immediately above,
one recovers the metric $g$. 
Namely, for $v,w \in T_mM$, choose $m^\prime \in \pi^{-1}(m)$ and
$v^\prime, w^\prime \in T_{m^\prime} OM$ such that $d\pi(v^\prime) = v$ and
$d\pi(w^\prime) = w$. Then $g(v,w) = \langle \tau(v^\prime), 
\tau(w^\prime) \rangle$.

Let $h : Spin(n) \rightarrow SO(n)$ be the double-covering homomorphism.
Giving a spin structure on $M$ means giving a principal $Spin(n)$-bundle $PM$
on $M$ such that $OM = PM \times_{Spin(n)} SO(n)$. The $1$-form $\tau$ lifts
to an $\R^n$-valued $1$-form $\tau^\prime$ on $PM$ with the properties that\\
1. If $V$ is a vertical vector field on $PM$ then $\tau^\prime(V) = 0$.\\
2. For all $\gamma \in Spin(n)$, 
$R_\gamma^* \tau^\prime = h(\gamma^{-1}) \cdot \tau^\prime$.\\
3. For all $f \in PM$, $\tau^\prime : T_fPM \rightarrow \R^n$ is onto. \\

Thus a Riemannian spin manifold consists of \\
1. The principal $Spin(n)$-manifold $PM$ on $M$ and\\
2. An $\R^n$-valued $1$-form $\tau^\prime$ on $PM$ satisfying
properties 1.-3. immediately above.\\
We can think of $PM$, as a topological fiber bundle, as being 
metric-independent.  Thus the notion of a spinor field on $M$ is also
metric-independent.  The metric only enters in defining the
$\R^n$-valued $1$-form $\tau^\prime$ on $PM$. In this way 
we can compare spinor fields
on two different Riemannian manifolds with the same underlying smooth 
structure. \\ \\
\noindent
{\bf Proof of Proposition \ref{prop1} :}

Suppose that the proposition is not true.  Then there is some $n \in \Z^+$,
some $\Lambda > 0$ and a sequence $\{\epsilon_i\}_{i = 1}^\infty$
of positive numbers such that\\
1. $\lim_{i \rightarrow \infty} \epsilon_i = 0$.\\
2. For each $i$, there is a connected closed $n$-dimensional spin manifold
$M_i$ with a Riemannian metric $g_i$ such that  
$-\epsilon_i \le K(M_i, g_i) \cdot \diam(M_i, g_i)^2 \le \Lambda$ and 
$\hA(M_i) \ne 0$.\\

By rescaling, we can assume that $\diam(M_i, g_i) = 1$. 
If $i$ is large enough then
$|K(M_i, g_i)| \le \Lambda$. We can write
\begin{equation} \label{Ahat}
\hA(M_i) = \int_{M_i} P(K(M_i, g_i)) \: d\vol(M_i)
\end{equation} 
for some explicit homogeneous polynomial $P$ in the
curvature tensor 
\cite[p. 231]{Lawson-Michelsohn (1989)}. Thus there is an explicit
number $v(n, \Lambda) > 0$
such that $\vol(M_i, g_i) 
\ge v(n, \Lambda)$, as otherwise we could conclude from
the integral formula that $\hA(M_i) = 0$. By Gromov's convergence theorem
and its elaborations 
\cite{Gromov-Lafontaine-Pansu (1978),Kasue (1989)}, there are \\
1. A smooth
manifold $M$ equipped with a metric $g_\infty$ which is $C^{1,\alpha}$-smooth
for all $0 < \alpha < 1$ and\\
2. A subsequence of $\{M_i\}_{i = 1}^\infty$, which we will relabel to again
call $\{M_i\}_{i = 1}^\infty$, and a sequence of diffeomorphisms
$F_i : M \rightarrow M_i$ such that $\lim_{i \rightarrow \infty}
F_i^* g_i = g_\infty$ in the $C^{1,\alpha}$-topology for all
$0 < \alpha < 1$.\\

Replacing $(M_i, g_i)$ by $(M, F_i^* g_i)$, we may assume that
the metrics $\{g_i\}_{i = 1}^\infty$ all live on the same
manifold $M$ (with $\hA(M) \ne 0$) 
and converge to $g_\infty$ in the $C^{1, \alpha}$-topology.
In particular, the Christoffel symbols of $g_\infty$ are locally 
$C^{0, \alpha}$ on $M$.
In fact, we may assume that for all $p \in [1, \infty)$, the sequence
$\{g_i\}_{i = 1}^\infty$ converges to $g_\infty$ in the Sobolev space
$L^{2,p}$ of covariant $2$-tensors on $M$
whose first two derivatives are $L^p$; a somewhat similar
case is treated in \cite[\S 2]{Anderson (1990)}.
Let $K_i$ denote the curvature tensor of $g_i$ and
let $K_\infty$ denote the curvature tensor of $g_\infty$, an $L^p$-tensor
for all $p \ge 1$. Then $\lim_{i \rightarrow \infty} K_i = K_\infty$ in
$L^p$ for all $p \ge 1$. In particular, $K_\infty \ge 0$ in the
sense of sectional curvatures.  Let $R_i$ denote the scalar curvature of 
$g_i$ and let $R_\infty \ge 0$ denote the scalar curvature of $g_\infty$.

As discussed above, we may take spinor fields to be sections of the
Hermitian vector bundle $S = PM \otimes_{Spin(n)} S_n$, regardless of
the Riemannian metric $g_i$.
Let $d\vol_i \in \Omega^n(M)$ be the volume form coming from $g_i$
and let $d\vol_\infty \in \Omega^n(M)$ be the volume form coming from
$g_\infty$.
Let $\nabla_i$ be the connection on $S$ coming from $g_i$ and
let $\nabla_\infty$ be the connection on $S$ coming from 
$g_\infty$. 
Then as $i \rightarrow \infty$, the tensor 
$\nabla_i - \nabla_\infty \in \End(S, S \otimes T^*M)$ converges
to zero in the $C^{0,\alpha}$-topology.
Let $D_i$ denote the Dirac operator on $S$ coming from $g_i$.
Let $H^0$ be the Hilbert space of $L^2$-spinors on $M$ with norm
\begin{equation}
\parallel \psi \parallel_{H^0}^2 = \int_M 
|\psi|^2 \: d\vol_\infty.
\end{equation}
Let $H^1$ be the Sobolev space of spinors on $M$ with norm
\begin{equation}
\parallel \psi \parallel_{H^1}^2 = \int_M \left( |\nabla_\infty \psi|^2 +
|\psi|^2 \right) \: d\vol_\infty.
\end{equation}

As $\hA(M) \ne 0$, the Atiyah-Singer index theorem implies that there
is a nonzero spinor field $\psi_i$ on $M$ such that $D_i \psi_i = 0$. 
We may assume that $\int_{M} |\psi_i|^2 \: d\vol_i = 1$.  From the 
Lichnerowicz formula,
\begin{equation} \label{lichno1}
0 = \int_{M} |D_i \psi_i|^2 \: d\vol_i = 
\int_M \left( |\nabla_i \psi_i|^2 + \frac{R_i}{4} \:
|\psi_i|^2 \right) \: d\vol_i.
\end{equation}
By our assumptions,
\begin{equation} \label{lichno2}
\int_M \frac{R_i}{4} \:
|\psi_i|^2 \: d\vol_i
\ge - \: \frac{n(n-1) \epsilon_i}{4}.
\end{equation}
Hence 
\begin{equation} \label{lichno3}
0 \le \int_M |\nabla_i \psi_i|^2 \: d\vol_i =
- \int_M \frac{R_i}{4} \:
|\psi_i|^2 \: d\vol_i
 \le \frac{n(n-1) \epsilon_i}{4}.
\end{equation}
Thus $\lim_{i \rightarrow \infty} \parallel \psi_i \parallel_{H^1} = 1$.
Taking a subsequence, we may assume that $\{ \psi_i \}_{\i=1}^\infty$
converges weakly to some $\psi_\infty \in H^1$. By compactness,
$\{ \psi_i \}_{\i=1}^\infty$ converges 
strongly to $\psi_\infty$ in $H^0$. Thus 
$\parallel \psi_\infty \parallel_{H^0} = 1$. Furthermore, for general
reasons,
\begin{equation}
\parallel \psi_\infty \parallel_{H^1} \: \le \: \lim_{i \rightarrow \infty}
\parallel \psi_i \parallel_{H^1} = 1.
\end{equation}
Hence
\begin{equation}
1 = \int_M |\psi_\infty|^2 \: d\vol_\infty \: \le \: 
\int_M \left( |\nabla_\infty \psi_\infty|^2 +
|\psi_\infty|^2 \right) \: d\vol_\infty \le 1.
\end{equation}
Thus $\nabla_\infty \psi_\infty = 0$. In particular, 
$| \psi_\infty |^2$ is a nonzero constant function on $M$.
Also, from (\ref{lichno3}), $\{\psi_i\}_{i=1}^\infty$ converges strongly
to $\psi_\infty$ in $H^1$.

As the $\widehat{A}$-genus is only nonzero in dimensions divisible by 
$4$, we may assume that $n > 2$. Then
$H^1$ embeds continuously in $L^{\frac{2n}{n-2}}$. Hence
$\lim_{i \rightarrow \infty} |\psi_i|^2 = |\psi_\infty|^2$ in
$L^{\frac{n}{n-2}}$. As $\lim_{i \rightarrow \infty} R_i = R_\infty$ in
$L^{\frac{n}{2}}$, (\ref{lichno3}) implies that 
\begin{equation}
\int_M \frac{R_\infty}{4} \:
|\psi_\infty|^2 \: d\vol_\infty = 
\lim_{i \rightarrow \infty} \int_M \frac{R_i}{4} \: |\psi_i|^2 \: d\vol_i = 0.
\end{equation} 
As $R_\infty \ge 0$, we conclude
that $R_\infty = 0$. Hence $K_\infty = 0$. As 
$\lim_{i \rightarrow \infty} K_i = K_\infty$ in
$L^{\frac{n}{2}}$, we obtain
\begin{equation}
\hA(M) = \lim_{i \rightarrow \infty} \int_{M} P(K_i) \: d\vol_i = 0.
\end{equation}
This is a contradiction. \qed \\ \\
{\bf Remark : }
The method of proof of Proposition \ref{prop1} gives the
following result.
\begin{proposition} \label{prop4}
For any $n \in \Z^+$ and any $\Lambda > 0$, there is an
$\epsilon(n, \Lambda) > 0$ such that if $M$ is a connected closed
$n$-dimensional Riemannian spin manifold with
$\Ric(M,g) \cdot \diam(M,g)^2 \ge - \epsilon(n, \Lambda)$ and
$|K(M,g)| \cdot \diam(M,g)^2 \le \Lambda$
then $\hA(M) = 0$ or $M$ admits a $C^{1,\alpha}$-metric
whose local holonomy group factorizes 
into products of $\{SU(m)\}_{m=2}^\infty$, 
$\{Sp(m)\}_{m=1}^\infty$, $Spin(7)$ and $G_2$. 
\end{proposition}

\begin{question} \label{quest3}
Given $n \in \Z^+$, is there an
$\epsilon(n) > 0$ such that if $M$ is a connected closed
$n$-dimensional Riemannian spin manifold with
$\Ric(M,g) \cdot \diam(M,g)^2 \ge - \epsilon(n)$
then $\hA(M) = 0$ or the frame bundle of $M$ admits a topological
reduction to a principal bundle whose local structure group
factorizes into products of $\{SU(m)\}_{m=2}^\infty$, 
$\{Sp(m)\}_{m=1}^\infty$, $Spin(7)$ and $G_2$? 
\end{question}

For example, $M = K3 \# (S^2 \times S^2)$ is a spin manifold with
$\hA(M) \ne 0$ but without an almost complex structure having 
$c_1 = 0$.  Does $M$ have almost-nonnegative-Ricci curvature? \\ \\
{\bf Remark :} One may think of trying to answer Question \ref{ques1} by
an extension of the Bochner method.  However, such an approach cannot
work, at least not directly.  For example, a flat torus is 
almost-nonnegatively-curved but, with the right spin structure, does have 
harmonic spinors.  It is just the index of its Dirac operator which vanishes.
Also, a nonflat nilmanifold has locally homogeneous metrics of 
constant negative scalar curvature, for
which the use of Lichnerowicz's formula is problematic.

\section{Proof of Proposition \ref{prop2}}
Part 1. of Proposition \ref{prop2} is proven in 
\cite[\S 2]{Fukaya-Yamaguchi (1992)}. More precisely,
if there is a metric $h$ on $B$ with $K(B, h) \cdot \diam(B, h)^2 > - \epsilon$
then there is a metric $g$ on $M$ with $K(M, g) \cdot \diam(M, g)^2 > 
- \epsilon$.

We now prove part 2. Put $b = \dim(B)$ and $n = \dim(N)$.
A fiber $N$ of $M$ has a tubular neighborhood which is diffeomorphic
to $N \times D^b$. As $M$ is spin, it follows that
$N \times D^b$ is spin and hence $N$ is spin. That is, $N$ is oriented and 
there
is a specific lifting of the oriented orthonormal frame bundle of $N$
to a principal $Spin(n)$-bundle $PN$. Let $G^\prime$ be the finite-index
subgroup of $G$
consisting of elements which are orientation-preserving on $N$ and 
lift to automorphisms of $PN$. Put $M^\prime = P \times_{G^\prime} N$.
Then $M^\prime$ is a finite covering of $M$. As the $\hA$-genus is
multiplicative under finite coverings, we may as well assume that
$G$ preserves the orientation and spin-structure of $N$. 

Suppose first that $N$ is not flat.  Choose a metric on $B$ and a 
connection on $P$. There is an induced metric on $M$.
As $N$ has nonnegative scalar curvature which is positive somewhere, 
it follows as in Lemma \ref{lemma1} that $\Ker(D_N)$ = 0.
If the fibers are made sufficiently small then $D_M$ is invertible
\cite[Proposition 4.41]{Bismut-Cheeger (1989)}.
The Atiyah-Singer index theorem then implies that $\hA(M) = \ind(D_M) = 0$.

Now suppose that $N$ is flat.  Then $N = T^k / F$ for some $k > 0$,
where $T^k$ has a flat metric, $F$ is a finite group of isometries
of $T^k$ and $T^k$ is a minimal such covering.  
Let $\rho : G \rightarrow \Isom(N)$ describe the action of
$G$ on $N$. Let $\Isom(T^k)^F$ denote the isometries of $T^k$ which 
commute with $F$. There is a homomorphism
$\Theta :  \Isom(T^k)^F \rightarrow \Isom(N)$. The
induced map on Lie algebras $\theta : \isom(T^k)^F \rightarrow \isom(N)$
is an isomorphism, as $\isom(N)$ is the Lie algebra of Killing vector
fields on $N$, each of which can be lifted to an $F$-invariant Killing 
vector field on $T^k$. As $\Ker(\Theta) = \cen(F)$, 
$\Theta$ restricts to an isomorphism between
$\Isom(T^k)^F_0$ and $\Isom(N)_0$, the connected components of the identity.

Put 
\begin{equation}
G^{\prime \prime} = \left\{ (g_1, g_2) \in \Isom(T^k)^F_0 \times G :
\Theta(g_1) = \rho(g_2) \right\}.
\end{equation}
There is a finite covering $P \times_{G^{\prime \prime}} 
T^k \rightarrow M$. As 
$\Isom(T^k)^F_0$ acts on $T^k$ by translations, it commutes with the action
of $T^k$ on itself by translations and so 
there is a nontrivial $T^k$-action on 
$P \times_{G^{\prime \prime}} T^k$.
By the Atiyah-Hirzebruch theorem \cite{Atiyah-Hirzebruch (1970)},
the $\hA$-genus of $P \times_{G^{\prime \prime}} T^k$ vanishes.
As the $\hA$-genus is multiplicative under finite covers, 
$\hA(M) = 0$. 

\section{Proof of Proposition \ref{prop3}}

Let $\overline{Z}$ be the union of the principal orbits for the action
of $G$ on $Z$. Put $\overline{B} = \overline{Z}/G$, a smooth manifold and put
$\overline{M} = \overline{Z} \times_G N$, a dense open subset of
$M$. There is a Riemannian submersion
$\pi : \overline{M} \rightarrow \overline{B}$ whose fiber over
$zG \in \overline{B}$ is $G_z \backslash N$.

To describe the geometry of $\overline{M}$ more explicitly, fix
$z \in \overline{Z}$. Let $N(G_z)$ denote the normalizer of $G_z$ in $G$.
Then $\overline{Z}$ is a fiber bundle over $\overline{B}$ with structure
group contained in $K = G_z \backslash N(G_z)$ 
\cite[Theorem 3.3]{Bredon (1972)}. That is, there is a principal $K$-bundle
$\overline{P}$ over $\overline{B}$ such that 
$\overline{Z} = \overline{P} \times_K (G_z\backslash G)$. 
Furthermore, $\overline{Z} \rightarrow \overline{B}$ is a Riemannian
submersion whose horizontal distribution comes from
a connection on $\overline{P}$. 
We note that although all of the $G$-orbits on $\overline{Z}$ are
diffeomorphic to $G_z\backslash G$, their Riemannian metrics
may vary from fiber to fiber. Topologically, we can write
$\overline{M} = \overline{P} \times_K (G_z \backslash N)$. 
The horizontal distribution
on the Riemannian submersion $\pi : \overline{M} \rightarrow \overline{B}$
again comes from the connection on $\overline{P}$. 
Metrically, the fibers of $\pi$ can be more
accurately written as $(G_z\backslash G) \times_G N$, with the orbit
$G_z\backslash G$ obtaining its metric from its embedding in 
$\overline{Z}$.\\ \\
{\bf Proof of 1.}
Let $g_Z$ and $g_N$ be the Riemannian metrics on $Z$ and $N$.
Let $K_0 > 0$ be such that $K(Z, g_Z) \ge -K_0$. 
For $j \ge 1$, consider the Riemannian metric 
$h_j = g_Z + j^{-2} g_N$ on $Z \times N$. Clearly
$K(Z \times N, h_j) \ge -K_0$ and $\diam(Z \times N, h_j) \le
\diam(Z, g_Z) + j^{-1} \diam(N, g_N)$. Let $(M, g_j) = (Z \times N, h_j)/G$ be 
the quotient metric on $M$. From the O'Neill formula 
\cite[Chapter 9]{Besse (1987)}, $K(M, g_j) \ge -K_0$. Clearly
$\diam(M, g_j) \le \diam(Z \times N, h_j)$. 

Let $\frac{N}{j}$ denote $(N, j^{-2} g_N)$. 
Let $\overline{g_j}$ denote
the restriction of $g_j$ to $\overline{M}$.
Then $(\overline{M}, \overline{g_j})$ is
obtained from the Riemannian submersion $\pi$ by 
changing the fiber from $(G_z \backslash G) \times_G N$
to $(G_z \backslash G) \times_G \frac{N}{j}$.
Let us concentrate on a given fiber $(G_z \backslash G) \times_G N$.
As $G_z \backslash N$ is a smooth manifold, there is a number $v_{min} > 0$
such that every $G_z$-orbit on $(N, g)$ has volume at least $v_{min}$.
Then 
\begin{equation}
\vol \left( (G_z \backslash G) \times_G N \right) \le
\frac{vol(G_z \backslash G) \cdot vol(N)}{vol(G_z \backslash G) \cdot v_{min}}
= \frac{vol(N)}{v_{min}}.
\end{equation}
Replacing $N$ by $\frac{N}{j}$ gives
\begin{equation}
\vol \left( (G_z \backslash G) \times_G \frac{N}{j} \right) \le
j^{- dim(G_z \backslash N)} \:
\frac{vol(N)}{v_{min}}.
\end{equation}
By assumption, $\dim(G_z \backslash N) > 0$. Thus
\begin{equation}
\lim_{j \rightarrow \infty} \vol \left(
(G_z \backslash G) \times_G \frac{N}{j} \right) = 0.
\end{equation}
It follows that
\begin{equation}
\lim_{j \rightarrow \infty} \vol(M, g_j) = 
\lim_{j \rightarrow \infty} \vol(\overline{M}, \overline{g_j})  = 0.
\end{equation}

We can think of the projection map $\pi : M \rightarrow Z/G$ as a singular
fibration whose fiber over $z^\prime G \in Z/G$ is $G_{z^\prime} \backslash N$.
From the same arguments as above, we see that $\lim_{j \rightarrow \infty}
(M, g_j) = Z/G$ in the Gromov-Hausdorff topology.  \\ \\
{\bf Proof of 2.} Without loss of generality, we may assume that
$\dim(M)$ is even.
As before, we fix $z \in \overline{Z}$.
Put $B = Z/G$,
\begin{equation}
Z^{sing} = \{ z^\prime \in Z : \dim(G_{z^\prime}) < \dim(G_{z})\},
\end{equation}
$B^{sing} = Z^{sing}/G$ and $M^{sing} = Z^{sing} \times_G N$.
Given $\epsilon > 0$, let 
$B^{sing}(\epsilon)$ be
the $\epsilon$-neighborhood of $B^{sing}$ in $B$, let $Z^{sing}(\epsilon)$
be its preimage in $Z$ and
put $M^{sing}(\epsilon) = Z^{sing}(\epsilon) \times_G N$. Put $M_1 = 
\overline{M^{sing}(\epsilon)}$, $M_2 = M - M^{sing}(\epsilon)$,
$W = \partial M_1 = \partial M_2$ and $B_2 = B - B^{sing}(\epsilon)$.
We note that $B_2$ is a smooth orbifold and that $M_2$ is a fiber bundle over
$B_2$.

By the O'Neill formula, $G_z\backslash N$
has a metric of nonnegative sectional curvature. Given this fact, it follows
from the Cheeger-Gromoll splitting theorem \cite{Cheeger-Gromoll (1971)} that
the condition that $G_z\backslash N$ be flat is topological in nature,
namely that $\pi_1(G_z\backslash N)$ have a free abelian subgroup of rank
$\dim(G_z\backslash N)$. We divide the proof of 2. into two cases.\\ \\
{\bf Case 1. $G_z\backslash N$ is not flat.}\\

We first prove a general result about the index of the Dirac operator on
a compact spin manifold-with-boundary. 
\begin{lemma} \label{vanishing}
Let $X$ be a compact even-dimensional
Riemannian spin manifold with boundary
$\partial X$. Let $D_X$ be the Dirac operator on $X$ with Atiyah-Patodi-Singer
boundary conditions 
\cite{Atiyah-Patodi-Singer (1975)}. Let $H_{\partial X} 
\in 
C^\infty(\partial X)$ be the mean curvature function.  (With our conventions,
if $X$ is the unit ball in $\R^n$, $n > 1$, then $H_{\partial X} > 0$.)
Suppose that $R_X \ge 0$ and $H_{\partial X} \ge 0$. Suppose that $R_X$ is
positive somewhere or $H_{\partial X}$ is positive somewhere. Then
$\ind(D_X) = 0$.
\end{lemma}
\begin{pf}
Let $\{e_j\}_{j=1}^n$ denote a local orthonormal frame on $X$. Let
$\gamma^j$ denote Clifford multiplication by $e_j$. With our conventions,
\begin{equation}
\gamma^i \gamma^j + \gamma^j \gamma^i = 2 \delta^{ij}.
\end{equation}
The Dirac operator on $X$ has the local form
\begin{equation}
D_X = -i \sum_{j=1}^n \gamma^j \: \nabla^X_{e_j}.
\end{equation}
Along $\partial X$, we take $e_n$ to be an inward-pointing unit normal
vector.  With respect to the decomposition $S = S^+ \oplus S^-$,
we can write
$\gamma^j = 
\begin{pmatrix}
0 & \sigma^j \\
\sigma^j & 0
\end{pmatrix}$ for $1 \le j \le n-1$ and
$\gamma^n = 
\begin{pmatrix}
0 & -i \\
i & 0
\end{pmatrix}$, where $\{\sigma_j\}_{j=1}^{n-1}$ are generators for
the Clifford algebra on $\R^{n-1}$. The Dirac operator on $\partial X$ has
the local form
\begin{equation}
D_{\partial X} = -i \sum_{j=1}^{n-1} \sigma^j \: \nabla^{\partial X}_{e_j}.
\end{equation}

Let $\psi$ be a spinor field on $X$. Let $\psi = \psi^+ + \psi^-$ be its
decomposition with respect to the $\Z_2$-grading on spinors.
Let $\psi_{\partial X}$ be its restriction to $\partial X$.  Let
$P^{\ge 0}$ be the projection onto the subspace of spinors on $\partial X$
spanned by eigenvectors of $D_{\partial X}$ of nonnegative eigenvalue, 
and similarly for
$P^{<0}$. The Atiyah-Patodi-Singer boundary conditions are
\begin{equation} \label{aps}
P^{\ge 0} \psi^+_{\partial X} = P^{< 0} \psi^-_{\partial X} = 0.
\end{equation}
These boundary conditions are usually considered when $X$ is a product near
the boundary, but 
one obtains an elliptic self-adjoint boundary condition for $D_X$ regardless
of whether or not $X$ is a product near the boundary.

Suppose that $D_X \psi = 0$. The Lichnerowicz equation
\begin{equation}
0 = D_X^2 \psi = \left(\nabla^X \right)^* \nabla^X \psi + \frac{R_X}{4} \: \psi
\end{equation}
is valid on the interior of $X$. Define a vector field $J = \sum_{j=1}^n
J^j e_j$ on $X$ by
$J^j = \langle \psi, \nabla^X_{e_j} \psi \rangle$. Then
\begin{align}
0 & = \int_X \langle \psi, \left(\nabla^X \right)^* \nabla^X \psi 
\rangle \: d\vol +
\int_X \frac{R_X}{4} \: |\psi|^2 \: d\vol \\
& = \int_X \langle \nabla^X \psi, \nabla^X \psi \rangle \: d\vol 
- \int_X \Div(J) \: d\vol
+ \int_X \frac{R_X}{4} \: |\psi|^2 \: d\vol \notag \\
& = \int_X \langle \nabla^X \psi, \nabla^X \psi \rangle \: d\vol
+ \int_X \frac{R_X}{4} \: |\psi|^2 \: d\vol  
+ \int_{\partial X} J^n \: d\vol \notag \\
\end{align} 
Now $D_X \psi = 0$ implies that $\nabla^X_{e_n} \psi = 
- \sum_{j=1}^{n-1} \gamma^n \gamma^j \: \nabla^X_{e_j} \psi$.
Hence
\begin{equation}
0  = \int_X \langle \nabla^X \psi, \nabla^X \psi \rangle \: d\vol
+ \int_X \frac{R_X}{4} \: |\psi|^2 \: d\vol  
- \int_{\partial X} \langle \psi, \sum_{j=1}^{n-1} \gamma^n \gamma^j \:
\nabla^X_{e_j} \psi \rangle \: d\vol.
\end{equation}
A computation gives that on $\partial X$,
\begin{equation}
\sum_{j=1}^{n-1} \gamma^n \gamma^j \: \nabla^X_{e_j} =
\sum_{j=1}^{n-1} \gamma^n \gamma^j \: \nabla^{\partial X}_{e_j} - 
(n-1) \: \frac{H_{\partial X}}{2}.
\end{equation}
Then
\begin{align} \label{mess}
0  = & \int_X \langle \nabla^X \psi, \nabla^X \psi \rangle \: d\vol
+ \int_X \frac{R_X}{4} \: |\psi|^2 \: d\vol \\ 
& - \int_{\partial X} \langle \psi, \sum_{j=1}^{n-1} \gamma^n \gamma^j \:
\nabla^{\partial X}_{e_j} \psi \rangle \: d\vol + 
 (n-1)  \int_{\partial X}
\frac{H_{\partial X}}{2} \: |\psi|^2 \: d\vol. \notag
\end{align}
The Atiyah-Patodi-Singer
boundary conditions (\ref{aps}) imply that 
\begin{equation}
- \int_{\partial X} \langle \psi, 
\sum_{j=1}^{n-1} \gamma^n \gamma^j \: \nabla^{\partial X}_{e_j} \psi \rangle 
\: d\vol \ge 0.
\end{equation}
As $R_X \ge 0$ and $H_{\partial X} \ge 0$, we obtain from
(\ref{mess}) that $\nabla^X \psi = 0$. In particular, $|\psi|^2$ is locally
constant on $X$. Equation (\ref{mess}), along with the fact that $R_X$ or
$H_{\partial X}$ is positive somewhere, implies that $|\psi|^2 = 0$.
Thus there are no nonzero solutions to $D_X \psi = 0$ and so $\ind(D_X) = 0$.
\end{pf}

Our strategy to prove the proposition in Case 1 is the following.  If $D_W$ is
invertible then the Atiyah-Patodi-Singer index theorem (and its
generalization to the case of nonproduct boundary) implies that
\begin{equation}
\ind(D_M) = \ind(D_{M_1}) + \ind(D_{M_2}).
\end{equation}
We will show that after shrinking the fiber metrics, we can apply
Lemma \ref{vanishing} to show that $\ind(D_{M_1}) = 0$. Then we will use
index theory techniques to show that $\ind(D_{M_2}) = 0$. 

\begin{lemma} \label{torpedo}
Define the metric $g_j$ on $M$ as in the proof of 
part 1. Then for large $j$, there is an $\epsilon > 0$ such that 
the submanifold
$M_1$ of $(M, g_j)$ has $R_{M_1} > 0$ and $H_{W} > 0$. 
\end{lemma}
\begin{pf}
This follows from computations as in \cite[Sect. 7-10]{Lewkowicz (1990)}.
We omit the details but give an illustrative example which has all of
the features of the general case.  Suppose that $z^\prime \in Z^{sing}$
is a $G$-fixed point.
By the equivariant tubular neighborhood theorem,
there is a $G$-vector space $V$ which is $G$-diffeomorphic to an
$\epsilon$-neighborhood of $z^\prime$ in $Z$
\cite[Theorem VI.2.2]{Bredon (1972)}. (In particular, for generic $v \in V$,
the $G$-stabilizer of $v$ is conjugate to $G_z$.) Then 
$U_M = V \times_{G} N$ is a neighborhood of
$\pi^{-1}(z^\prime G) \subset M$. Consider the case when
$N$ is a homogeneous space $G/H$.

If $N$ has positive scalar curvature then it is easy to see that from the
O'Neill formula that for large $j$, $(U_M, g_j)$ has positive scalar
curvature.  Suppose, to take the other extreme, that $N$ is flat. Taking
a finite cover, we may assume that $N = T^k$ and that $G$ acts
on $T^k$ through a homomorphism $\rho^\prime : G \rightarrow
T^k$.  As $U_M$ is
a smooth manifold, $\Ker(\rho^\prime)$ must act trivially on $V$.
Put $\widetilde{G} = G/\Ker(\rho^\prime)$. There are
induced homomorphisms $\widetilde{\rho} : \widetilde{G} \rightarrow \Aut(V)$
and $\widetilde{\rho}^\prime : \widetilde{G} \rightarrow
T^k$, with $\widetilde{\rho}^\prime$ being an inclusion. Then 
$U_M = V \times_{\widetilde{G}} T^k$. Taking a finite cover,
we may assume that $\widetilde{G}$ is connected. If $V$ is flat 
then one can check that for large $j$, $(U_M, g_j)$ 
has a positive scalar 
curvature function whose value at $[0,t] \in V \times_{\widetilde{G}} T^k$ is
$O(j^2)$. This can be seen intuitively by the fact that
$\widetilde{\rho}$ reduces into trivial $\R$-factors and at least
one nontrivial $\R^2$-factor.  If $V = \R^2$ and ${\widetilde{G}} = T^k = S^1$
then $\R^2 \times_{S^1} S^1$ has a torpedo
shape which becomes more curved at the tip as the $S^1$-factor shrinks.

In the general case, the torpedo effect ensures that if
$j$ is large enough and $\epsilon$ is small enough then
$(U_M, g_j)$ will have positive
scalar curvature.  In fact, for large $j$ we can take 
$\epsilon = O(j^{-(\frac{1}{2} + \alpha)})$ for any $\alpha > 0$.
As $B^{sing} \subset B$ has codimension at least two,
the mean curvature of $\partial M_1$
is positive for large $j$.  Doing a similar procedure
for a finite collection of $z^\prime \in Z^{sing}$, we can deal with all of the
strata of $Z^{sing}$. The lemma follows.
\end{pf}

As $G_z \backslash N$ is not flat, each fiber of the fiber bundle
$M_2 \rightarrow B_2$ has a nonnegative scalar curvature
function which is positive somewhere.  Then by the Lichnerowicz formula,
the Dirac operator on each fiber is invertible. 
For large $j$, the geometry of a fiber
$(G_z \backslash G) \times_G \frac{N}{j}$ 
is asymptotically that of $\frac{G_z \backslash N}{j}$.
If $j$ is large and $0 < \alpha < \frac{1}{2}$ then it follows as in
\cite[Proposition 4.41]{Bismut-Cheeger (1989)} that
$D_W$ is also invertible. By Lemmas \ref{vanishing} and \ref{torpedo},
$\ind(D_{M_1}) = 0$.

We now show that $\ind(D_{M_2}) = 0$.
Let $[0, \delta) \times W \subset M_2$ be a neighborhood of $W$
such that if $u \in [0, \delta)$ is the coordinate function then
$\partial_u$ is a unit length vector field whose flow generates unit-speed
geodesics which are normal to $W$.
We can write the metric near $W$ as $du^2 + h(u)$, where $h(u)$ is a metric
on $W$. Let $F : [0, \infty) \rightarrow [0,1]$ be a smooth nondecreasing
function such that $F(x)$ is identically zero for $x$ near zero and
identically one if $x \ge \frac{1}{2}$. For $v \in [0,1]$, define
$f_v : [0, \infty) \rightarrow [0, \infty)$ by
\begin{equation}
f_v(u) =  
\begin{cases}
u & \text{if $v = 0$,} \\
u \: F( \frac{u}{\delta v} ) & \text{if $v \in (0,1]$.} 
\end{cases}
\end{equation}
Let $M_2(v)$ be the manifold $M_2$ with the metric
$du^2 + h(f_v(u))$ on $[0, \delta) \times W$. Then $M_2(0)$ is the same
as $M_2$ with the original metric and $M_2(1)$ has a product metric near
its boundary. For all $v \in [0,1]$, $\partial M_2(v)$ is isometric to $W$.
Then the Dirac
operators on $M_2(v)$, with Atiyah-Patodi-Singer boundary conditions,
form a continuous family of Fredholm operators and so have constant index
with respect to $v$. Thus for computational purposes, we may assume that
$M_2$ is a product near the boundary.

By the Atiyah-Patodi-Singer index theorem \cite{Atiyah-Patodi-Singer (1975)},
\begin{equation} \label{eqstart}
\ind(D_{M_2}) = \int_{M_2} \widehat{A}(\nabla^{TM_2}) - \frac{1}{2} \:
\eta_W.
\end{equation}
From \cite[Theorems 4.35 and 4.95]{Bismut-Cheeger (1989)}, we have an equality
in $\Omega^*(B_2)$ :
\begin{equation}
d\widetilde{\eta}_{M_2} = \int_{G_z \backslash N} 
\hA \left(\nabla^{T(G_z \backslash N)} \right).
\end{equation}
(Strictly speaking, we have to generalize
the results of \cite{Bismut-Cheeger (1989)} from smooth fiber bundles to
fiber bundles with orbifold base.  Such a generalization is straightforward.
We will not give the details here.)

Also,
\begin{align}
\lim_{j \rightarrow \infty} \int_{M_2} \widehat{A}(\nabla^{TM_2}) & =
\int_{B_2} \int_{G_z \backslash N} \hA \left(\nabla^{T(G_z \backslash N)} 
\right) \wedge \hA \left(\nabla^{TB_2} \right) =
\int_{B_2} d\widetilde{\eta}_{M_2} \wedge \hA \left(\nabla^{TB_2} \right) \\
& = \int_{\partial B_2} \widetilde{\eta}_{W} \wedge \hA \left(
\nabla^{T\partial B_2} \right). \notag
\end{align}
On the other hand, by 
\cite[Theorems 4.35 and 4.95]{Bismut-Cheeger (1989)},
\begin{equation} \label{eqend}
\lim_{j \rightarrow \infty} \frac{1}{2} \: \eta_W =
\int_{\partial B_2} \widetilde{\eta}_{W} \wedge \hA \left( 
\nabla^{T\partial B_2}
\right).
\end{equation}
Combining equations (\ref{eqstart})-(\ref{eqend}) gives that
$\ind(D_{M_2}) = 0$. This proves the proposition in Case 1.\\ \\
{\bf Case 2. $G_z \backslash N$ is flat.} \\

We will show that there is a finite cover of $M$ with a nontrivial
$S^1$-action.  It then follows from the Atiyah-Hirzebruch theorem
\cite{Atiyah-Hirzebruch (1970)} that $\hA(M) = 0$.

The various fibers of the fiber bundle $\overline{M} \rightarrow \overline{B}$
are all flat.  They are not necessarily all isometric.  However, they are
all affine-equivalent.

Write $G_z \backslash N = F \backslash T^k$, 
where $k > 0$, $F \subset \Aff(T^k)$ is a finite
group of affine diffeomorphisms of $T^k$ and $T^k$ is a minimal such covering.
Let $\rho : K \rightarrow
\Aff(G_z \backslash N)$ describe the action of $K$ on $G_z \backslash N$.
Let $\Aff(T^k)^F$ denote the affine diffeomorphisms of $T^k$ which commute
with $F$ and let $\Aff(T^k)^F_0$ denote the connected component of the
identity.  There is a homomorphism 
$\Theta : \Aff(T^k)^F \rightarrow \Aff(G_z \backslash N)$. Put
\begin{equation}
K^{\prime} = \left\{ (k_1, k_2) \in \Aff(T^k)^F_0 \times K :
\Theta(k_1) = \rho(k_2) \right\}.
\end{equation}
Put
$\overline{M}^\prime = \overline{P} \times_{K^{\prime}} T^k$, a finite
cover of $\overline{M}$. Then $\overline{M}^\prime$ is a fiber bundle
over $\overline{B}$ with fiber $K \times_{K^\prime} T^k$, a finite disjoint
union of tori. There is a nontrivial $T^k$-action on $\overline{M}^\prime$.
Similarly, we want to show that there is a finite cover 
$M^{\prime}$ of $M$ with a nontrivial $T^k$-action.  The problem is that
$\overline{M}^\prime$ may not extend to a finite cover of $M$. However,
we will show that the disjoint union of a certain number of copies of it
does extend.

Choose $z^\prime \in Z - \overline{Z}$. Let $G_{z^\prime} \subseteq G$ be its
stabilizer subgroup. Then $G_z \subset G_{z^\prime}$ and
$G_{z^\prime} \backslash N$ is a smooth manifold.  
There is a Riemannian submersion
$G_z \backslash N \rightarrow G_{z^\prime} \backslash N$ with fiber
$G_z \backslash G_{z^\prime}$. As $G_z \backslash N$ is flat and
both $G_{z^\prime} \backslash N$ and 
$G_z \backslash G_{z^\prime}$ are nonnegatively curved, one sees from the
homotopy groups that $G_{z^\prime} \backslash N$ and
$G_z \backslash G_{z^\prime}$ must be flat.  As $G_z \backslash G_{z^\prime}$
is a globally homogeneous space, it must be a disjoint union of tori
of dimension $\dim(G_{z^\prime}) - \dim(G_z)$.

By the equivariant tubular neighborhood theorem, 
there is a finite-dimensional
real vector space $V$, a
representation $\rho : G_{z^\prime} \rightarrow \Aut(V)$ and
a neighborhood $U_Z$ of the $G$-orbit of $z^\prime$ such that
$U_Z$ is 
$G$-diffeomorphic to $V \times_{G_{z^\prime}} G$. Then
$U_M = U_Z \times_G N = V \times_{G_{z^\prime}} N$ is a neighborhood of
$\pi^{-1}(z^\prime G) \subset M$.

As $G_{z^\prime} \backslash N$ is flat, we can write it as
$F^\prime \backslash T^{k^\prime}$, where $F^\prime$ is a finite group
of affine diffeomorphisms of $T^{k^\prime}$ and 
$T^{k^\prime}$ is a minimal such covering. Let $s : T^{k^\prime} \rightarrow
F^\prime \backslash T^{k^\prime}$ be the projection map. 
Consider the fiber bundle $F \backslash T^k \stackrel{r}{\rightarrow}
F^\prime \backslash T^{k^\prime}$ with fiber 
$G_z \backslash G_{z^\prime}$. Put 
\begin{equation}
C = \{ (t, t^\prime) \in (F \backslash T^k) \times T^{k^\prime} :
r(t) = s(t^\prime) \}.
\end{equation}  
Equivalently, $C = r^* T^{k^\prime}$, as shown in the diagram
\begin{equation} \label{diagram1}
\begin{array}{ccc}
C & \rightarrow & T^{k^\prime} \\
\downarrow & & \downarrow s \\
F \backslash T^k & \stackrel{r}{\rightarrow} & F^\prime \backslash 
T^{k^\prime}.
\end{array}
\end{equation}
We claim that $C$ is a disjoint union of $k$-dimensional tori.
To see this,  put $\Gamma = \pi_1(G_z \backslash N)$ and 
$\Gamma^\prime = \pi_1(G_{z^\prime} \backslash N)$. Then
we have a diagram of exact sequences :
\begin{equation} \label{diagram2}
\begin{array}{lllllllllll}
 &  &  &  & 
&  & 1 &  & 
 &  &  \\
 &  &  &  &
&  & \downarrow &  & 
 &  &  \\
 &  &  &  &
& & \Z^{k^\prime} & & 
 &  &  \\
 &  &  &  &
&  & \downarrow \beta &  & 
 &  &  \\
1 & \rightarrow & \pi_1(G_z \backslash G_{z^\prime}) & \rightarrow & \Gamma
& \stackrel{\alpha}{\rightarrow} & \Gamma^\prime & \rightarrow & 
\pi_0(G_z \backslash G_{z^\prime}) & \rightarrow & 1. \\
 &  &  &  &
&  & \downarrow &  & 
 &  &  \\
 &  &  &  &
& & F^\prime &  & 
 &  &  \\
 &  &  &  &
&  & \downarrow &  & 
 &  &  \\
 &  &  &  &
& & 1 &  & 
 &  &
\end{array}
\end{equation}
Now
\begin{equation}
\pi_1(C) = \{ (\gamma, \zeta) \in \Gamma \times \Z^{k^\prime} : 
\alpha(\gamma) = \beta(\zeta)\}. 
\end{equation}
Projecting $\pi_1(C)$ on $\Gamma$ or $\Z^{k^\prime}$, it 
follows from (\ref{diagram2}) that there are exact sequences
\begin{equation} \label{es1}
1 \rightarrow \pi_1(C) \rightarrow \Gamma
\end{equation}
and 
\begin{equation} \label{es2}
1 \rightarrow \pi_1(G_z \backslash G_{z^\prime}) \rightarrow \pi_1(C)
\rightarrow \Z^{k^\prime}.
\end{equation}
Put $A = \Image(\pi_1(C) \rightarrow
\Z^{k^\prime})$, a free abelian subgroup of finite index.
Then (\ref{es2}) is equivalent to the short exact sequence
\begin{equation} \label{es3}
1 \rightarrow \pi_1(G_z \backslash G_{z^\prime}) \rightarrow \pi_1(C)
\rightarrow A \rightarrow 1.
\end{equation}
From (\ref{es1}), $\pi_1(C)$ is isomorphic to a subgroup of $\Gamma$.
One can see that it is of finite index in $\Gamma$ and so has
polynomial growth of degree $k$. This implies that the sequence
(\ref{es3}) splits and that $\pi_1(C) =  \pi_1(G_z \backslash G_{z^\prime})
\times A \cong \Z^k$. As $C$ is flat, it must be the
disjoint union of $m$ copies of $T^k$ for some $m > 0$.

Let $\delta : N \rightarrow G_{z^\prime} \backslash N$ be the projection
map.  Put
\begin{equation}
R = \{ (n, t^\prime) \in N \times T^{k^\prime} : \delta(n) = s(t^\prime)\}. 
\end{equation} 
That is, $R = \delta^* T^{k^\prime}$,
as shown in the diagram
\begin{equation} \label{diagram3}
\begin{array}{ccc}
R & \rightarrow & T^{k^\prime} \\
\downarrow & & \downarrow s\\
N & \stackrel{\delta}{\rightarrow} & G_{z^\prime} \backslash N,
\end{array}
\end{equation}
whose vertical arrows are finite coverings.
There is an action of $G_{z^\prime}$ on $R$ by
$g (n, t^\prime) = (gn, t^\prime)$, with $G_z \backslash R = C$
(compare (\ref{diagram1})).

Let $[K : K^\prime]$ denote the index of $K^\prime$ in $K$.
Let $[K : K^\prime] \: R$ denote the disjoint union of $[K : K^\prime]$ copies
of $R$.
Consider the finite covering 
$V \times_{G_{z^\prime}} [K : K^\prime] \: R
\rightarrow U_M$.
We claim that
this extends the covering $m \overline{M}^\prime \rightarrow \overline{M}$
over $U_M$. To see this,
note that we have a diagram of fiber bundles
\begin{equation} \label{diagram4}
\begin{array}{ccccc}
G_z \backslash G_{z^\prime} & \rightarrow & 
V \times_{G_{z}} [K : K^\prime] \: R & \rightarrow &
V \times_{G_{z^\prime}} [K : K^\prime] \: R \\
& & \downarrow & & \downarrow \\
G_z \backslash G_{z^\prime} & \rightarrow & 
V \times_{G_{z}} N & \rightarrow &
V \times_{G_{z^\prime}} N.
\end{array}
\end{equation}
Let $\overline{V}$ be the set of points in $V$ whose stabilizer group is
conjugate to $G_z$, a dense open subset of $V$. Put $\overline{U}_M =
\overline{M} \cap U_M =
\overline{V} \times_{G_{z^\prime}} N$, a dense open subset of $U_M$. 
Let $\overline{U}^\prime_M$ be the pre-image of $\overline{U}_M$ under
the covering $\overline{M}^\prime \rightarrow \overline{M}$.
Then
(\ref{diagram4}) restricts to
\begin{equation} \label{diagram5}
\begin{array}{ccccc}
G_z \backslash G_{z^\prime} & \rightarrow & 
\overline{V} \times_{G_{z}} [K : K^\prime] \: R & \rightarrow &
\overline{V} \times_{G_{z^\prime}} [K : K^\prime] \: R \\
& & \downarrow & & \downarrow \\
G_z \backslash G_{z^\prime} & \rightarrow & 
\overline{V} \times_{G_{z}} N & \rightarrow &
\overline{V} \times_{G_{z^\prime}} N.
\end{array}
\end{equation}
That is, we have a diagram of fiber bundles
\begin{equation} \label{diagram6}
\begin{array}{ccccc}
& & m [K : K^\prime] \: T^k & & m [K : K^\prime] \: T^k \\
& & \downarrow & & \downarrow \\
G_z \backslash G_{z^\prime} & \rightarrow & 
\overline{V} \times [K : K^\prime] \: C & \rightarrow &
\overline{V} \times_{G_{z^\prime}} [K : K^\prime] \: R \\
& & \downarrow & & \downarrow \\
G_z \backslash G_{z^\prime} & \rightarrow & 
\overline{V} \times ({G_{z}} \backslash N) & \rightarrow &
\overline{U}_M.
\end{array}
\end{equation}
By the constructions, it follows that the right-hand-column of
(\ref{diagram6}) is the same as
\begin{equation}
\begin{array}{c}
m K \times_{K^\prime} T^k \\
\downarrow \\
m \overline{U}^\prime_M \\
\downarrow \\
\overline{U}_M.
\end{array}
\end{equation}
Thus $V \times_{G_{z^\prime}} [K : K^\prime] \: R
\rightarrow U_M$ does
extend the covering $m \overline{M}^\prime \rightarrow \overline{M}$
over $U_M$. Furthermore, the obvious $T^k$-action on $m \overline{U}^\prime_M$
comes from the $T^k$-action on $\overline{V} \times [K : K^\prime] \: C$,
which extends to the $T^k$-action on $V \times_{G_z} [K : K^\prime] \: R
= V \times [K : K^\prime] \: C$,
which pushes down to a $T^k$-action on 
$V \times_{G_{z^\prime}} [K : K^\prime] \: R$. Of course, the $T^k$-action
on $V \times_{G_{z^\prime}} [K : K^\prime] \: R$ may not be free.

Repeating the process for a finite number of $z^\prime$'s whose $G$-orbits
exhaust the singular
$G$-strata of $Z$, we end up with a finite covering $M^\prime
\rightarrow M$. The preimage of $\overline{M}$ in $M^\prime$ is the disjoint
union of a finite number of copies of $\overline{M}^\prime$ and so has
a nontrivial $T^k$-action.  From the nature of the above extension procedure,
we know that it extends to a $T^k$-action on $M^\prime$.
The proposition follows.


\begin{thebibliography}{9}

\bibitem{Anderson (1990)} M. Anderson, ``Convergence and Rigidity of
Manifolds under Ricci Curvature Bounds'', Inv. Math. 102, p. 429-445
(1990)


\bibitem{Atiyah-Hirzebruch (1970)}
M. Atiyah and F. Hirzebruch, ``Spin-Manifolds and Group Actions'',
in \underline{Essays on Topology and Related Topics}, Springer, New York,
p. 18-28 (1970) 

\bibitem{Atiyah-Patodi-Singer (1975)} 
M. Atiyah, V. Patodi and I. Singer, ``Spectral Asymmetry and
Riemannian Geometry I'', Math. Proc. of Camb. Phil. Soc. 77, p. 43-69 (1975)

\bibitem{Besse (1987)} A. Besse, \underline{Einstein Manifolds},
Springer-Verlag, New York (1987)

\bibitem{Bismut-Cheeger (1989)} 
J.-M. Bismut and J. Cheeger, ``$\eta$-Invariants and Their
Adiabatic Limits'', J. Amer. Math. Soc. 2, p. 33-70 (1989)

\bibitem{Bredon (1972)} G. Bredon, \underline{Introduction to Compact
Transformation Groups}, Academic Press, New York (1972)

\bibitem{Cheeger-Gromoll (1971)} J. Cheeger and D. Gromoll,
``The Splitting Theorem for Manifolds of Nonnegative Ricci Curvature'',
J. Diff. Geom. 6, p. 119-128 (1971)

\bibitem{Fukaya-Yamaguchi (1992)} K. Fukaya and T. Yamaguchi,
``The Fundamental Groups of Almost Nonnegatively-Curved Manifolds'',
Ann. Math. 136, p. 253-333 (1992)

\bibitem{Gallot (1983)} S. Gallot, ``In\'egalit\'es Isop\'erim\'etriques,
Courbure de Ricci et Invariants G\'eom\'etriques II'', C. R. Acad. Sci.
S\'er. I Math. 296, p. 365-368 (1983)

\bibitem{Gromov (1981)} M. Gromov, ``Curvature, Diameter and Betti Numbers'',
Comm. Math. Helv. 56, p. 179-195 (1981)

\bibitem{Gromov-Lafontaine-Pansu (1978)} M. Gromov, J. Lafontaine and
P. Pansu, \underline{Structures M\'etriques pour les Vari\'et\'es
Riemanniennes}, CEDIC, Paris (1978)

\bibitem{Grove-Petersen-Wu (1990)} K. Grove, P. Petersen V and J.-Y. Wu,
``Geometric Finiteness Theorems via Controlled Topology'', Inv. Math.
99, p. 205-213 (1990)

\bibitem{Kasue (1989)} A. Kasue, ``A Convergence Theorem for Riemannian
Manifolds and Some Applications'', Nagoya Math. J. 114, p. 21-51 (1989) 

\bibitem{Lawson-Michelsohn (1989)} H.B. Lawson and M.-L. Michelsohn,
\underline{Spin Geometry}, Princeton University Press, Princeton, New Jersey
(1989)

\bibitem{Lewkowicz (1990)} M. Lewkowicz, ``Positive Scalar Curvature and
Local Actions of Nonabelian Lie Groups'', J. Diff. Geom. 31, p. 29-45 (1990)
 
\bibitem{Petersen-Wei-Ye (1995)} P. Petersen V, G. Wei and R. Ye,
``Controlled Geometry via Smoothing'', preprint,
http://xxx.lanl.gov/abs/dg-ga/9508012 (1995)

\bibitem{Weiss (1996)} M. Weiss, ``Curvature and Finite Domination'',
Proc. Amer. Math. Soc. 124, p. 615-622 (1996)
\end{thebibliography}
\end{document}